\documentclass[10pt]{article}
\usepackage[utf8]{inputenc}
\usepackage{latexsym,amssymb,amsmath,url,authblk}

\def\N{\mathbb{N}}

\def\C{\mathbb{C}}

\def\proof{\par\noindent{\em Proof. }}
\def\eproof{\hfill{$\Box$}\bigskip}

\def\ds{\dots}

\newtheorem{thm}{Theorem}[section]

\title{Arnol'd's limit and the Lagrange inversion}
\author{Martin Klazar\footnote{{\tt klazar@kam.mff.cuni.cz}}}
\date{\today}

\begin{document}

\maketitle

\begin{abstract}
We show how to prove by means of the Lagrange inversion the limit of Arnol'd 
that 
$$
\lim_{x\to0}\frac{\sin(\tan x)-\tan(\sin x)}{\arcsin(\arctan x)-\arctan(\arcsin x)}=1\,.
$$
In fact, we obtain a~more general result in terms of formal power series.
\end{abstract}

\section{Arnol'd's limit}\label{sec_limi}

In \cite[p.~21]{arno}, see \cite{arno_eng} for English translation, Vladimir~I.~Arnol'd (1937--2010) poses to the reader, in his sarcastic style, the following challenge.
\begin{quote}
Here is an example of a~problem that would be solved by people like Barrow, Newton or Huygens in a~couple 
of  minutes, but contemporary mathematicians, in my opinion, are unable to solve it quickly (in any case I have not yet seen a~mathematician who could cope with it quickly): evaluate
$$
\lim_{x\to0}\frac{\sin(\tan x)-\tan(\sin x)}{\arcsin(\arctan x)-\arctan(\arcsin x)}\;.
$$
\end{quote}
On the one hand, Arnol'd is right, in a~sense. At least this author could not do with the limit anything simple
for a~long time. In one July night 2025 he only 
laboriously computed, when writing the textbook/monograph \cite{klaz}, that the numerator and denominator are both
$$
{\textstyle
-\frac{1}{30}x^7+o(x^7)\ \ (x\to0)\,.
}
$$
On the other hand, K.~Honn pointed out in recent preprint \cite{honn} that Arnol'd's geometric resolution of the limit, as given in \cite[note~8 on p.~85]{arno}, is problematic, and corrected it. Later the same night it suddenly occurred to this author: 
analytic functions, hence power series, and inverse functions, how does it rhyme?
The Lagrange inversion formula (LIF)! Then he (almost) cursed, because this should have 
occurred to him at least 25 years ago, and rushed to the notebook (well, later 
in the morning). In the next section we show how  LIF can be used to establish 
a~result much more general than Arnol'd's limit.

\section{The Lagrange inversion formula}

For $n\in\N_0$ ($=\{0,1,\ds\}$) and a~formal power series 
$$
{\textstyle
f=f(x)=\sum_{n\ge0}a_nx^n\ \ (\in\C[[x]])
}
$$
we write $[x^n]f$ for the coefficient $a_n$ of $x^n$. We say that $f(x)$
in $\C[[x]]$ is {\em regular} if $[x^0]f=0$ and $[x^1]f\ne0$. It is well 
known that every regular $f(x)$ has a~unique {\em inverse}, a~regular formal power series $f^{\langle-1\rangle}(x)$ such that
$$
f(f^{\langle-1\rangle}(x))=f^{\langle-1\rangle}(f(x))=x\,.
$$
LIF, in the form we apply here, says the following.

\begin{thm}[LIF]\label{thm_LIF}
If $f(x)$ is a~regular formal power series and $n\in\N$ ($=\{1,2,\ds\}$), then
$$
[x^n]f^{\langle-1\rangle}=n^{-1}\cdot[x^{n-1}]\Big(\frac{x}{f(x)}\Big)^n\,.
$$
\end{thm}
For more information on LIF see \cite{gess}.

For a~nonzero formal power series $f(x)$ we write $\mathrm{ord}(f)$ to denote the minimum
$n\in\N_0$ such that $[x^n]f\ne0$, and we
set
$$
m(f)=\big[x^{\mathrm{ord}(f)}\big]f\ \ (\in\C^*)\,.
$$
The following theorem is our main result.

\begin{thm}\label{thm_main}
Suppose that $f(x)$ and $g(x)$ are two distinct regular formal power series such that $[x^1]f=[x^1]g=1$. Then
$$
\mathrm{ord}(f-g)=\mathrm{ord}(g^{\langle-1\rangle}-f^{\langle-1\rangle})\,\text{ and }\,
m(f-g)=m(g^{\langle-1\rangle}-f^{\langle-1\rangle})\ \ (\in\C^*)\,.
$$
\end{thm}
\proof
If $f(x)=\sum_{n\ge0}a_nx^n$ and $g(x)=\sum_{n\ge0}b_nx^n$, then 
$a_0=b_0=0$, $a_1=b_1=1$,
$$
\mathrm{ord}(f-g)=k\,\text{ and }\,
m(f-g)=a_k-b_k\ne0\,,
$$
where $k\in\N_0$, $k\ge2$, is uniquely determined by the condition that $a_k\ne b_k$ but 
$a_n=b_n$ for every $n\in\N_0$ with $n<k$.
Using Theorem~\ref{thm_LIF} we get for every $n\in\N$ that
\begin{eqnarray*}
[x^n]f^{\langle-1\rangle}&=&{\textstyle
\frac{1}{n}[x^{n-1}]\big(\frac{x}{f(x)}\big)^n}\\
&=&{\textstyle
\frac{1}{n}[x^{n-1}]
\Big(\sum_{j=0}^{n-1}(-1)^j\big(
a_2x+a_3x^2+\ds+a_nx^{n-1}\big)^j\Big)^n\,,
}
\end{eqnarray*}
and we have the same formula for $[x^n]g^{\langle-1\rangle}$, only $a_2$, $a_3$, $\ds$, $a_n$ are replaced with $b_2$, $b_3$, $\ds$, $b_n$, respectively. 
If $n<k$, it follows that
$$
[x^n]f^{\langle-1\rangle}=[x^n]g^{\langle-1\rangle}\,.
$$
For $n=k$ we get from the above formula that
\begin{eqnarray*}
[x^k]f^{\langle-1\rangle}&=&{\textstyle
\frac{1}{k}[x^{k-1}]
\Big(1-\big(a_2x+a_3x^2+\ds+a_kx^{k-1}\big)^1\,+}\\
&&+\,\mathrm{poly}(a_2x,\,a_3x^2,\,\ds,\,a_{k-1}x^{k-2})\Big)^k={\textstyle
\frac{1}{k}\cdot k\cdot 1^{k-1}\cdot (-a_k)^1}\,+\\
&&+\,
\mathrm{poly}_0(a_2,\,a_3,\,\ds,\,a_{k-1})\\
&=&-a_k+\mathrm{poly}_0(a_2,\,a_3,\,\ds,\,a_{k-1})\,,
\end{eqnarray*}
where $\mathrm{poly}(\cdots)$ and 
$\mathrm{poly}_0(\cdots)$
are rational polynomials in the stated variables and $\mathrm{poly}(\cdots)$ has 
zero constant term. Again, we have the same formula, with the same polynomials $\mathrm{poly}(\cdots)$ and 
$\mathrm{poly}_0(\cdots)$, for $[x^k]g^{\langle-1\rangle}$, only $a_2$, $a_3$, $\ds$, $a_k$ are replaced with $b_2$, $b_3$, $\ds$, $b_k$, respectively. Thus
$$
\mathrm{ord}(g^{\langle-1\rangle}-f^{\langle-1\rangle})=k\,\text{ and \,}m(g^{\langle-1\rangle}-f^{\langle-1\rangle})=a_k-b_k\,.
$$
The theorem is proven.
\eproof

\noindent
In particular, if $f(x)$, respectively $g(x)$, is the Taylor series of 
$\sin(\tan x)$, respectively $\tan(\sin x)$, with center $0$ (which can be 
interpreted as obtained by composing two formal power series), Theorem~\ref{thm_main} gives that 
Arnol'd's limit equals $1$. Why is $f(x)\ne g(x)$? This is because, for example, 
$$
\text{$\lim_{x\to\pi/2}\sin(\tan x)$ does not exist but $\lim_{x\to\pi/2}\tan(\sin x)=\tan(1)$}\,.
$$

\section{Concluding remarks}

We hope to return to Arnol'd's limit in near future in the next versions of this note. For 
example, a~natural idea is to adapt the 
Arnol'd--Honn geometric proof to functions that are not analytic and 
satisfy only some smoothness conditions.

\medskip\noindent
{\em Department of Applied Mathematics\\
Faculty of Mathematics and Physics\\
Charles University\\
Malostransk\'e n\'am\v es\'\i\ 25\\
118 00 Praha\\
Czechia}

\end{document}